\documentclass{amsart}

\usepackage{ae}
\usepackage{amssymb}

\usepackage[alphabetic]{amsrefs}
\usepackage{mathrsfs}

\newcommand{\sconn}{{\operatorname{sc}}}
\newcommand{\Kalg}{K^{\operatorname{alg}}}
\newcommand{\Ksep}{K^{\operatorname{sep}}}

\newcommand{\C}{\mathcal{C}}
\newcommand{\e}{\varepsilon}

\newcommand{\lie}[1]{\mathfrak{#1}}
\newcommand{\glie}{\lie{g}}
\newcommand{\hlie}{\lie{h}}
\newcommand{\clie}{\lie{c}}
\newcommand{\sllie}{\lie{sl}}
\newcommand{\mlie}{\lie{m}}
\newcommand{\slie}{\lie{s}}

\newcommand{\plie}{\lie{p}}

\newcommand{\Lie}{\operatorname{Lie}}
\newcommand{\G}{\mathbf{G}}

\newcommand{\F}{\mathbf{F}}
\newcommand{\Gm}{{\G_m}}
\newcommand{\Ga}{{\G_a}}

\newcommand{\Proj}{\mathbf{P}}
\newcommand{\Stab}{\operatorname{Stab}}

\newcommand{\Ad}{\operatorname{Ad}}

\newcommand{\Int}{\operatorname{Int}}

\newcommand{\red}{{\operatorname{red}}}
\newcommand{\image}{\operatorname{im}}
\newcommand{\congruent}{\equiv}

\newcommand{\gcr}{$G$-cr}
\newcommand{\T}{\mathcal{T}}

\newcommand{\Z}{\mathbf{Z}}

\newcommand{\SL}{\operatorname{SL}}
\newcommand{\GL}{\operatorname{GL}}
\newcommand{\SP}{\operatorname{Sp}}
\newcommand{\SO}{\operatorname{SO}}

\newcommand{\tensor}{\otimes}

\newtheorem{prop}{Proposition}

\newtheorem{lem}[prop]{Lemma}
\newtheorem{theorem}[prop]{Theorem}
\newtheorem{cor}[prop]{Corollary} 

\theoremstyle{remark}
\newtheorem{example}[prop]{Example}

\newtheorem{rem}[prop]{Remark}

\numberwithin{equation}{subsection}

\begin{document}
\title{Completely reducible $\operatorname{SL}(2)$-homomorphisms}
\author{George J. McNinch}
\address{Department of Mathematics,
         Tufts University,
         503 Boston Avenue,
         Medford, MA 02155,
         USA}
\email{george.mcninch@tufts.edu}
\author{Donna M. Testerman}
\address{Institut de g\'eom\'etrie, alg\`ebre et topologie,
  B\^atiment BCH,
  \'Ecole Polytechnique F\'ed\'erale de Lausanne,
  CH-1015 Lausanne,
  Switzerland}

\email{donna.testerman@epfl.ch} \thanks{Research of McNinch supported
  in part by the US National Science Foundation through DMS-0437482.}
\thanks{Research of Testerman supported in part by the Swiss National
  Science Foundation grant PP002-68710.}

\date{October 18, 2005}

\begin{abstract}
  Let $K$ be any field, and let $G$ be a semisimple group over $K$.
  Suppose the characteristic of $K$ is positive and is very good for
  $G$.  We describe all group scheme homomorphisms $\phi:\SL_2 \to G$
  whose image is geometrically $G$-completely reducible -- or \gcr\ --
  in the sense of Serre; the description resembles that of irreducible
  modules given by Steinberg's tensor product theorem. In case $K$ is
  algebraically closed and $G$ is simple, the result proved here was
  previously obtained by Liebeck and Seitz using different methods.  A
  recent result shows the Lie algebra of the image of $\phi$ to be
  geometrically \gcr; this plays an important role in our proof.
\end{abstract}

\maketitle

\section{Introduction}
\label{sec:intro}

Let $K$ be an arbitrary field of characteristic $p>0$. By a scheme we
mean a separated $K$-scheme of finite type. An algebraic group will
mean a smooth and affine $K$-group scheme; a subgroup will mean a
$K$-subgroup scheme, and a homomorphism will mean a $K$-homomorphism.
A smooth group scheme $G$ is said to be reductive if $G_{/\Kalg}$ is
reductive in the usual sense -- i.e.  it has trivial unipotent radical --
where $\Kalg$ is an algebraic closure of $K$.  The Lie algebra $\glie
= \Lie(G)$ may be regarded as a scheme over $K$; we permit ourselves
to write $\glie$ for the set of $K$-points $\glie(K)$.

For $G$ a reductive group, a subgroup $H \subset G$ is said to be
geometrically $G$-completely reducible -- or \gcr\ -- if whenever $k$
is an algebraically closed field containing $K$ and $H_{/k}$ is
contained in a parabolic $k$-subgroup $P$ of $G_{/k}$, then $H_{/k} \subset
L$ for some Levi $k$-subgroup $L$ of $P$; see \S\ref{sub:Lie-cr} for
more details. The notion of \gcr\ was introduced by J-P. Serre; see
e.g. \cite{serre-sem-bourb} for more on this notion. It is our goal
here to describe all homomorphisms $\phi:\SL_2 \to G$ whose image is
geometrically \gcr; this we achieve under some assumptions on $G$
which are described in \S \ref{sub:strongly-standard}. For the
purposes of this introduction, let us suppose that $G$ is semisimple.
Then our assumption is: the characteristic of $K$ is \emph{very good}
for $G$ (again see \S \ref{sub:strongly-standard} for the precise definition of
a very good prime).

Let $F:\SL_2 \to \SL_2$ be the Frobenius endomorphism obtained by base
change from the Frobenius endomorphism of $\SL_{2/\F_p}$; cf. \S
\ref{sub:frob} below.  We say that a collection of homomorphisms
$\phi_0,\phi_1,\dots,\phi_r:\SL_2 \to G$ is
\emph{commuting} if
\begin{equation*}
  \image \phi_i \subset C_G(\image \phi_j) \quad \text{for all $0 \le i \ne j \le r$}.
\end{equation*}
Let $\vec \phi = (\phi_0,\dots,\phi_r)$ where the $\phi_i$ are
commuting homomorphisms $\SL_2 \to G$, and let $\vec n = (n_0 < \dots
<n_r)$ where the $n_i$ are non-negative integers. Then the data $(\vec
\phi,\vec n)$ determines a homomorphism $\Phi_{\vec \phi,\vec n}:\SL_2
\to G$ given for every commutative $K$-algebra $\Lambda$ and every $g
\in \SL_2(\Lambda)$ by the rule
\begin{equation*} 
    g \mapsto \phi_0(F^{n_0}(g)) \cdot \phi_1(F^{n_1}(g)) \cdots \phi_r(F^{n_r}(g)).
\end{equation*}
We say that $\Phi=\Phi_{\vec \phi,\vec n}$ is the
\emph{twisted-product homomorphism} determined by $(\vec \phi,\vec
n)$.

A notion of \emph{optimal} homomorphisms $\SL_2 \to G$ was introduced
in \cite{mcninch-optimal}; see \S \ref{sub:optimal} for the precise
definition. When $G$ is a $K$-form of $\GL(V)$ or $\SL(V)$, a
homomorphism $f:\SL_2 \to G$ is optimal just in case the
representation $(f_{/\Ksep},V)$ is restricted and semisimple, where
$\Ksep$ is a separable closure of $K$; see Remark
\ref{rem:SL(V)-optimal}.  We will say that the list of commuting
homomorphisms $\vec\phi = (\phi_0,\phi_1,\dots,\phi_r)$ is
\emph{optimal} if each $\phi_i$ is an optimal homomorphism.

\begin{theorem}
  \label{theorem:main}
  Let $G$ be a semisimple group for which the characteristic is very good, and let
  $\Phi:\SL_2 \to G$ be a homomorphism.  If the image of $\Phi$ is
  geometrically \gcr, then there are commuting optimal homomorphisms
  $\vec \phi = (\phi_0,\dots,\phi_r)$ and non-negative integers $\vec
  n = (n_0 < n_1 < \cdots < n_d)$ such that $\Phi$ is the
  twisted-product homomorphism determined by $(\vec \phi,\vec n)$.
  Moreover, $\vec \phi$ and $\vec n$ are uniquely determined by
  $\Phi$.
\end{theorem}

We actually prove the theorem for the strongly standard reductive
groups described below in \ref{sub:strongly-standard}; see Theorem
\ref{theorem:real-main}.

In case $K$ is algebraically closed and $G$ is simple, this theorem
was obtained by Liebeck and Seitz \cite{liebeck-seitz}*{Theorem 1};
cf. Remark \ref{rem:good-vs-optimal} to see that the notion of
restricted -- or good -- $A_1$-subgroup used in \cite{liebeck-seitz}
is ``the same'' as the notion of optimal homomorphism used here.

Note that Liebeck and Seitz prove a version of Theorem
\ref{theorem:main} where $\SL_2$ is replaced by any quasisimple group
$H$.  If $G$ is a split classical group over $K$ in good
characteristic, the more general form of Theorem \ref{theorem:main}
found in \cite{liebeck-seitz} is a consequence of Steinberg's tensor
product theorem \cite{JRAG}*{Cor. II.3.17}; cf.
\cite{liebeck-seitz}*{Lemma 4.1}.  The proof given by Liebeck and
Seitz of Theorem \ref{theorem:main} for a quasisimple group $G$ of
exceptional type relies instead on detailed knowledge of the subgroup
structure -- in particular, of the maximal subgroups -- of $G$; see
e.g.  \cite{liebeck-seitz}*{Theorem 2.1, Proposition 2.2, and \S 4.1}
for the case $H=\SL_2$. In contrast, when $p>2$, our proof uses in an
essential way the complete reducibility of the Lie algebra of a \gcr\
subgroup of $G$ \cite{mcninch-lie-cr}; cf. the proofs of Lemma
\ref{lem:slie-not-in-P}, Proposition
\ref{prop:lambda0-lambda-centralizes-slie}, and  Lemma
\ref{lem:parabolic-inclusion} [when $p=2$, we have essentially just
used the proof of Liebeck and Seitz].

We obtain also the converse to Theorem \ref{theorem:main}, though we
do so only under a restriction on $p$.  Write $h(G)$ for the maximum
value of the Coxeter number of a simple $k$-quotient of $G_{/k}$,
where $k$ is an algebraically closed field containing $K$.

\begin{theorem}
  \label{theorem:converse}
  Let $G$ be semisimple in very good characteristic, and  suppose that $p >
  2h(G) -2$, let $\vec \phi = (\phi_0,\dots,\phi_d)$ be commuting
  optimal homomorphisms $\SL_2 \to G$, and let $\vec n = (n_0 < n_1 <
  \cdots < n_d)$ be non-negative integers. Then the image of the
  twisted-product homomorphism $\Phi:\SL_2 \to G$ determined by $(\vec
  \phi,\vec n)$ is geometrically \gcr.
\end{theorem}

Again, this result is proved for a more general class of reductive
groups; see Theorem \ref{theorem:real-converse}.

The assumption on $p$ made in the last theorem is unnecessary if $G$
is a classical group -- or a group of type $G_2$ -- in good
characteristic; see Remark \ref{rem:converse-classical}. However, it
is not clear to the authors how to eliminate the prime restriction in
general.

The first named author would like to acknowledge the hospitality of
the Centre Interfacultaire Bernoulli at the \'Ecole Polytechnique
F\'ed\'erale de Lausanne during a visit in June 2005; this visit
permitted much of the collaboration which led to the present
manuscript.

\section{Preliminaries}
\label{sec:prelim}

\subsection{Reduced subgroups}
\label{sub:reduced}
Let $k$ be a perfect field -- in the application we take $k$ to be
algebraically closed. Let $B$ be a group scheme of finite type over $k$.

\begin{lem}
  \label{lem:reduced}
  There is a unique smooth subgroup $B_\red \subset B$ which has the same
  underlying topological space as $B$. If $A$ is any smooth group
  scheme over $k$ and $f:A \to B$ is a $k$-morphism, then $f$ factors
  in a unique way into a $k$-morphism $A \to B_\red$ followed by the
  inclusion $B_\red \to B$.
\end{lem}

\begin{proof}
  Use \cite{liu}*{Prop. 2.4.2} to find the reduced $k$-scheme $B_\red$
  with the same underlying topological space as $B$; the result just
  quoted then yields the uniqueness of $B_\red$.  It is clear that
  $B_\red$ is a $k$-group scheme, and the assertion about $A$ and $f$
  follows from \emph{loc. cit.} Prop 2.4.2(d).  Since $k$ is perfect,
  apply \cite{KMRT}*{Prop. 21.9} to see that a $k$-group is smooth if
  and only if it is geometrically reduced if and only if it is
  reduced. Thus $B_\red$ is indeed smooth .
\end{proof}

We are going to consider later some group schemes which we do not
\emph{a priori} know to be smooth, and we want to choose maximal tori
in these group schemes. The following example explains why in those
cases we first extend scalars to an algebraically closed field (see e.g.
\S \ref{sub:geom} below).

\begin{example}
  If $B$ is a group scheme over an imperfect field $K$, and if $k$ is
  a perfect field containing $K$, then a maximal torus of $B_{/k,\red}$
  need not arise by base-change from a $K$-subgroup of $B$. Let us
  give an example.

  Let $A = \Gm \ltimes \Ga$ where $\Gm$ acts on $\Ga$ ``with weight
  one''; i.e.  $K[A] = K[T,T^{-1},U]$ where the comultiplication
  $\mu^*$ is given by
  \begin{equation*}
    \mu^*(T^{\pm 1}) = (T \tensor T)^{\pm 1} 
    \quad \text{and} \quad 
    \mu^*(U) = U \tensor T + 1 \tensor U.    
  \end{equation*}

  Suppose that $K$ is not perfect, and let $L = K(\beta)$ where
  $\beta^p = \alpha \in K$ but $\beta \not \in K$. Consider the subgroup
  scheme $B \subset A$ defined by the ideal 
  $I = (\alpha T^p - U^p - \alpha) \lhd K[A]$.

  If $k$ is a perfect field containing $K$, notice that the image
  $\bar{f} \in k[B]$ of $f=\beta T - U - \beta \in k[A]$ satisfies
  $\bar{f}^p = 0$ but $\bar{f} \ne 0$; thus $B_{/k}$ is not reduced.
  The subgroup $B_{/k,\red} \subset A_{/k}$ is defined by $J = (\beta T - U
  - \beta)$, so that $B_{/k,\red} \simeq \Gm_{/k}$ is a torus. The
  group of $k$-points $B_{/k,\red}(k) \subset  A(k)$ may be described
  as:
  \begin{equation*}
    B_{/k,\red}(k) = \{(t,\beta t - \beta) \in \Gm(k) \ltimes \Ga(k) \mid t \in k^\times\}.
  \end{equation*}
  Note that $B_{/k,\red}$ does not arise by base change from a
  $K$-subgroup of $A$, e.g. since the intersection $B_{/k,\red}(k)
  \cap A(K)$ consists only in the identity element [where the
  intersection takes place in the group $A(k)$].
\end{example}


\subsection{Cocharacters and parabolic subgroups}
A \emph{cocharacter} of an algebraic group $A$ is a homomorphism
$\gamma:\Gm \to A$. We write $X_*(A)$ for the set of cocharacters of
$A$.

A linear representation $(\rho,V)$ of $A$ yields a linear
representation $(\rho \circ \gamma,V)$ of $\Gm$ which in turn is
determined by the morphism
\begin{equation*}
  (\rho\circ \gamma)^*:V \to K[\Gm] \tensor_K V = K[t,t^{-1}] \tensor_K V.
\end{equation*}
Then $V$ is the direct sum of the weight spaces
\begin{equation}
  \label{eq:weight-spaces}
  V(\gamma;i) = \{v \in V \mid (\rho \circ \gamma)^*v = t^i \tensor v \}
\end{equation}
for $i \in \Z$. 

Consider now the  reductive group $G$.  If $\gamma \in X_*(G)$,
then
\begin{equation*}
  P_G(\gamma)= P(\gamma) = 
  \{ x \in G \mid \lim_{t \to 0} \gamma(t)x\gamma(t^{-1})\text{\ exists}\}
\end{equation*}
is a parabolic subgroup of $G$ whose Lie algebra is $\plie(\gamma) =
\sum_{i \ge 0} \glie(\gamma;i)$; see e.g. \cite{springer-LAG}*{\S 3.2} for
the notion of limit used here.  Moreover, each parabolic subgroup of
$G$ has the form $P(\gamma)$ for some cocharacter $\gamma$; for all
this cf. \cite{springer-LAG}*{3.2.15 and 8.4.5}.

We note that $\gamma$ ``exhibits'' a Levi decomposition of
$P=P(\gamma)$. Indeed, $P(\gamma)$ is the semi-direct product
$C_G(\gamma) \cdot U(\gamma)$, where $U(\gamma) = \{x \in P \mid
\lim_{t\to 0} \gamma(t)x\gamma(t^{-1}) = 1\}$ is the unipotent radical of
$P(\gamma)$, and the reductive subgroup $C_G(\gamma) =
C_G(\gamma(\G_m))$ is a Levi factor in $P(\gamma)$; cf.
\cite{springer-LAG}*{13.4.2}.

\subsection{Complete reducibility, Lie algebras}
\label{sub:Lie-cr}

Let $G$ be a reductive group, and write $\glie$ for its Lie algebra.

A smooth subgroup $H \subset G$ is geometrically \gcr\ if whenever $k$
is an algebraically closed field containing $K$ and $H_{/k} \subset P$
for a parabolic $k$-subgroup $P \subset G_{/k}$, then $H_{/k} \subset
L$ for some Levi $k$-subgroup $L \subset P$.

Similarly, if $\hlie \subset \glie$ is a Lie subalgebra, we say that
$\hlie$ is geometrically \gcr\ if whenever $k$ is an algebraically
closed field containing $K$ and $P \subset G_{/k}$ is a parabolic
$k$-subgroup with $\hlie_{/k} = \hlie \tensor_K k \subset \Lie(P)$,
then $\hlie_{/k} \subset \Lie(L)$ for some Levi $k$-subgroup $L
\subset P$.

 \begin{lem}
   Let $X$ and $Y$ be schemes of finite type over $K$, and let $f:X
   \to Y$ be a ($K$-) morphism. The following are equivalent:
   \begin{enumerate}
   \item[i)] $f$ is surjective,
     \item[ii)] $f_{/k}:X(k) \to Y(k)$ is surjective for all algebraically closed fields
       $k$ containing $K$, and
     \item[iii)] $f_{/k}:X(k) \to Y(k)$ is surjective for some algebraically closed field
       $k$ containing $K$.
   \end{enumerate}
 \end{lem}

 \begin{proof}
   This follows from \cite{DG}*{I\ \S 3.6.10}
 \end{proof}

 \begin{lem}
   \label{lem:indep-of-acf}
   Fix an algebraically closed field $k$ containing $K$.  Let $G$ be a
   reductive group, let $J \subset G$ be a smooth subgroup, and let
   $\hlie \subset \glie$ be a Lie subalgebra.  Then
   \begin{enumerate}
   \item $J$ is geometrically $G$-cr if and only if $J_{/k}$ is $G_{/k}$-cr.
   \item $\hlie$ is geometrically $G$-cr if and only if $\hlie_{/k}$ is $G_{/k}$-cr.
   \end{enumerate}
 \end{lem}

 \begin{proof}
   We prove (1); the proof of (2) is essentially the same.  We are
   going to apply the previous Lemma.

   First let $\mathcal{P}$ be the scheme of all parabolic subgroups of
   $G$, and let $Y = \mathcal{P}^J$ be the fixed point scheme for the
   action of $J$; thus $Y$ is the closed subscheme of those parabolic subgroups
   containing $J$.  \footnote{For assertion (2), one should instead
     regard $\mathcal{P}$ as the scheme of parabolic subalgebras of
     $\glie$, which may be regarded as a closed subscheme of a product
     of Grassman schemes $\operatorname{Gr}_d(\glie)$ for various $d$.
     Now the subscheme $X \subset \mathcal{P}$ of parabolic
     subalgebras containing $\hlie$ is the intersection of
     $\mathcal{P}$ with the subscheme $Z$ of the product of Grassman
     schemes consisting of those subspaces containing $\hlie$. Since
     $Z$ is closed in the product, $Y$ is closed in $\mathcal{P}$. Similar
     remarks apply to the definition of the subscheme $Y \subset
     \mathcal{PL}$ to be given in the next paragraph.}

   Let also $\mathcal{PL}$ be the scheme such that for each
   commutative $K$-algebra $\Lambda$, the $\Lambda$-points
   $\mathcal{PL}(\Lambda)$ are the pairs $P \supset L$ where $P$ is a
   parabolic of $G_{/\Lambda}$ and $L$ is a Levi subgroup of $P$; cf.
   \cite{SGA3}*{Exp. XXVI, \S 3.15}.  Let $X = (\mathcal{PL})^J$ be the
   scheme of those pairs $P \supset L$ where $L$ contains $J$.

   There is an evident morphism $\mathcal{PL} \to \mathcal{P}$ given
   by $(P \supset L) \mapsto P$; cf.  \cite{SGA3}*{Exp. XXVI, \S
     3.15}.  By restriction one gets a morphism $f:X \to Y$. Then $f$
   is surjective if and only if $J$ is \gcr, and (1) follows from the
   preceding Lemma.
\end{proof}

 \begin{prop}
   \label{prop:gcr-levi}
   Let $G$ be reductive, and let $M \subset G$ be a Levi subgroup.
   Suppose that $J \subset M$ is a smooth subgroup, and that $\hlie
   \subset \Lie(M)$ is a Lie subalgebra.  Then $J$ is geometrically
   \gcr\ if and only if $J$ is geometrically $M$-cr and $\hlie$ is
   geometrically \gcr\ if and only if $\hlie$ is geometrically $M$-cr.
 \end{prop}

 \begin{proof}
   For the proof, it is enough to suppose that $K$ is algebraically
   closed.
   The proof for $J$ is found in \cite{bmr}*{Theorem 3.10}.   
   The proof for $\hlie$ is deduced from
   \cite{serre-sem-bourb}*{2.1.8}; see \cite{mcninch-lie-cr}*{Lemma
     2} for the argument.
 \end{proof}

The following theorem was proved in \cite{mcninch-lie-cr}.
\begin{theorem}
  \label{theorem:lie-cr}
  Let $H \subset G$ be a smooth subgroup which is  geometrically \gcr.
  Then $\hlie = \Lie(H)$ is  geometrically \gcr.
\end{theorem}

We recall a similar result of B. Martin \cite{serre-sem-bourb}*{Th\'eor\`eme 3.6}.
\begin{theorem}[Martin]
  \label{theorem:martin}
  Let $H \subset G$ be a smooth subgroup which is geometrically \gcr, and
  let $H' \lhd H$ be a smooth normal subgroup. Then $H'$ is
  geometrically \gcr\ as well.
\end{theorem}

 Finally, we note:
 \begin{lem}
   \label{lem:gcr-isogeny}
   Let $\pi:G \to G_1$ be a central isogeny where $G_1$ is a second
   reductive group, let $J \subset G$ be a smooth subgroup, and let
   $\hlie \subset \glie = \Lie(G)$ be a Lie subalgebra. Then
   \begin{enumerate}
   \item $J$ is geometrically \gcr\ if and only if $\pi(J)$ is
     geometrically $G_1$-cr, and
   \item $\hlie$ is geometrically \gcr\ if and only if $d\pi(\hlie)$
     is geometrically $G_1$-cr.
   \end{enumerate}
 \end{lem}

 \begin{proof}
   We may and will suppose that $K$ is algebraically closed for the proof.
   It is clear that $J$ is contained in a parabolic subgroup
   $P$ of $G$ if and only if $\pi(J)$ is contained in the parabolic
   subgroup $\pi(P)$ of $G_1$, and similarly $\hlie$ is contained in
   $\Lie(P)$ if and only if $d\pi(\hlie)$ is contained in
   $d\pi(\Lie(P)) = \Lie(\pi(P))$, the result follows since $P \mapsto
   \pi(P)$ determines a bijection between the parabolic subgroups of
   $G$ and those of $G_1$.
 \end{proof}

\subsection{Strongly standard reductive groups}
\label{sub:strongly-standard}

If $G$ is geometrically quasisimple with absolute root system $R$
\footnote{The absolute root system of $G$ is the root system of
  $G_{/\Ksep}$ where $\Ksep$ is a separable closure of $K$.}, the
characteristic $p$ of $K$ is said to be a bad prime for $R$ in the
following circumstances: $p=2$ is bad whenever $R \not = A_r$, $p=3$
is bad if $R = G_2,F_4,E_r$, and $p=5$ is bad if $R=E_8$.  Otherwise,
$p$ is good.  [Here is a more intrinsic definition of good prime: $p$
is good just in case it divides no coefficient of the highest root in
$R$].

If $p$ is good, then $p$ is said to be very good provided that either
$R$ is not of type $A_r$, or that $R=A_r$ and $r \not
\congruent -1 \pmod p$.

There is a possibly inseparable central isogeny \footnote{Indeed, the
  center of the reductive group $G$ is a smooth subgroup scheme; this
  follows e.g.  from \cite{SGA3}*{II Exp.  XII Th\'eor\`eme 4.1} since
  for reductive $G$, the center is the same as the ``centre
  r\'eductif''. The radical $R(G)$ is the maximal torus of the center
  of $G$, so $R(G)$ is a smooth torus, and we take $T=R(G)$ in
  \eqref{eq:cover-G}.  Now, multiplication gives a central isogeny $G'
  \times R(G) \to G$ where $G'$ is the derived group of $G$. So
  \eqref{eq:cover-G} follows from the corresponding result for
  semisimple groups; see e.g.  \cite{KMRT}*{Theorems 26.7 and 26.8} or
  \cite{tits-weiss}*{Appendix (42.2.7)}.}
\begin{equation}
  \label{eq:cover-G}
   \quad \prod_{i=1}^r G_i \times T \to G
\end{equation}
for some torus $T$ and some $r \ge 1$, where for $1 \le i \le r$ there
is an isomorphism $G_i \simeq R_{L_i/K} H_i$ for a finite
separable field extension $L_i/K$ and a geometrically simple, simply
connected $L_i$-group scheme $H_i$; here, $R_{L_i/K} H_i$
denotes the ``Weil restriction'' of $H_i$ to $K$.

Then $p$ is good, respectively very good, for $G$ if and only if that
is so for $H_i$ for every $1 \le i \le r$. Since the $H_i$ are
uniquely determined by $G$ up to central isogeny, the notions of good
and very good primes depend only on the central isogeny class of the
derived group $(G,G)$. Moreover, these notions are geometric in the
sense that they depend only on the group $G_{/k}$ for an algebraically closed field
$k$ containing $K$.

One says that a smooth $K$-group $D$ is of multiplicative type if
$D_{/K'}$ is diagonalizable for some algebraic extension $K \subset
K'$; i.e. that $D_{/K'} \simeq \operatorname{Diag}(\Gamma)$ for some
commutative group $\Gamma$. See \cite{JRAG}*{I.2.5} for the definition
of $\operatorname{Diag}(\Gamma)$ -- it is implicitly defined in
\cite{springer-LAG}*{Corollary 3.2.4} as well. A torus is of multiplicative
type, as is any finite, smooth, commutative subgroup all of whose geometric
points are semisimple.

If $G$ is reductive and if $D \subset G $ is a
subgroup of multiplicative type, then $C_G(D)$ is a reductive subgroup
containing a maximal torus of $G$ -- use \cite{SGA3}*{II Exp. XI, Cor
  5.3} to see that $C_G(D)$ is smooth, use
\cite{springer-LAG}*{Theorem 3.2.3} to see that $D_{K'}$ lies in a maximal
torus of $G_{K'}$, and finally use \cite{springer-steinberg}*{II. \S 4.1}
to see that $C_G(D)$ is reductive.

Consider reductive groups of the form
\begin{equation*}
  (*) \quad H=H_1 \times S
\end{equation*}
where $H_1$ is a semisimple group for which the characteristic of $K$
is very good, and where $S$ is a torus. We say that $G$ is
\emph{strongly standard} if there is a group $H$ as in $(*)$, a
subgroup of multiplicative type $D\subset H$, and a separable isogeny between
$G$ and the reductive subgroup $C_H(D)$ of $H$.

\begin{rem}
\label{rem:strongly-standard-remark}
This definition of \emph{strongly standard} is more general than that
given e.g. in \cite{mcninch-optimal}. It follows from Proposition
\ref{prop:strongly-standard-properties} below that the main result of
\emph{loc. cit.} in fact applies to a strongly standard group in this
stronger sense.
\end{rem}

\begin{prop}
  \label{prop:strongly-standard-properties}
  Let $G$ be strongly standard.
  \begin{enumerate}
  \item If $D \subset G$ is a subgroup of multiplicative type, then the
    reductive group $C_G(D)$ is strongly standard.
  \item The characteristic of $K$ is good for the derived group of
    $G$, and there is a non-degenerate, $G$-invariant bilinear form on
    $\Lie(G)$.
  \item Each conjugacy class and each
  adjoint orbit is separable. In particular, if $g \in G(K)$ and $X
  \in \glie(K)$, then $C_G(g)$ and $C_G(X)$ are smooth.
  \footnote{In older language, these centralizers are defined over $K$.}
  \end{enumerate}
\end{prop}

\begin{rem}
  The centralizers considered in (3) -- and elsewhere in this paper --
  are the \emph{scheme-theoretic} centralizers. Thus e.g. $C_G(X)$ is
  the group scheme with $\Lambda$-points $C_G(X)(\Lambda) = \{g \in
  G(\Lambda) \mid \Ad(g)X = X\}$ for each commutative $K$-algebra
  $\Lambda$.
\end{rem}

\begin{proof}[Proof of Proposition \ref{prop:strongly-standard-properties}]
  For the proofs of (1) and (2), we may replace $G$ by a separably
  isogenous group and suppose $G$ to be the centralizer of a subgroup of
  multiplicative type $D_1 \subset H$ where $H$ has the form $(*)$.

  For (1), note that since $D_1$ centralizes $D$, the group $D_2 = D
  \cdot D_1 \subset H$ is of multiplicative type and (1) is immediate.

  To prove (2), note first that the characteristic of $K$ is good for
  the derived group of $G$ by \cite{mcninch-optimal}*{Lemma 1}(2).
  Now, there is a non-degenerate $H$-invariant bilinear form $\beta$
  on a group $H$ of the form $(*)$ by \cite{mcninch-optimal}*{Lemma
    1}(1). Moreover, it suffices to see that the restriction of
  $\beta$ to $\Lie(G)$ is nondegenerate after making a field extension; thus, we may
  suppose that $D_1 \simeq \operatorname{Diag}(\Gamma)$.  We have
  $\Lie(H) = \bigoplus_{\gamma \in \Gamma} \Lie(H)_\gamma$ where $D_1$
  acts on $\Lie(H)_\gamma$ through $\gamma$; see e.g.  \cite{JRAG}*{\S
    I.2.11}.  The subspaces $\Lie(H)_\gamma$ and $\Lie(H)_\tau$ are
  evidently orthogonal unless $\gamma \cdot \tau = 1$ in $\Gamma$
  \footnote{We are writing $\Gamma$ multiplicatively}.  Since $\Lie(G)
  = \Lie(H)^{D_1} = \Lie(H)_1$, the restriction of $\beta$ to
  $\Lie(G)$ must remain non-degenerate.

  In view of (2), the proof of \cite{mcninch-optimal}*{Prop. 5} yields
  (3).
\end{proof}

\subsection{Nilpotent elements and associated cocharacters}
\label{sub:nilpotent}

Let $G$ be a reductive group, and let $X \in \glie=\Lie(G)$ be
nilpotent. A cocharacter $\Psi \in X_*(G)$ is said to be associated
with $X$ if the following conditions hold:

\begin{enumerate}
\item[(A1)] $X \in \glie(\Psi;2)$, and
\item[(A2)] there is a maximal torus $S$ of $C_G(X)$ such that
  $\Psi \in X_*(L_1)$ where $L=C_G(S)$ and $L_1 = (L,L)$ is its derived group.
\end{enumerate}

Assume now that $G$ is strongly standard.

\begin{prop}
  \label{prop:assoc-cochar}
  Let $X \in \glie$ be nilpotent. 
  \begin{enumerate}
  \item There is a cocharacter $\Psi$ associated with $X$.
  \item If $\Psi$ is associated to $X$ and $P=P(\Psi)$ is the
    parabolic determined by $\Psi$, then $C_G(X) \subset P$.  In
    particular, $\clie_\glie(X) \subset \Lie(P)$.
  \item If $\Psi,\Phi \in X_*(G)$  are associated
    with $X$, then $\Psi = \Int(x) \circ \Phi$ for some $x \in C_G(X)(K)$.
  \item The parabolic subgroups $P(\Psi)$ for cocharacters $\Psi$
    associated with $X$ all coincide.
  \end{enumerate}
\end{prop}

\begin{proof}
  (1) is \cite{mcninch-rat}*{Theorem 26}, (2) is
  \cite{jantzen-nil}*{Prop. 5.9}. (3) follows from
  \cite{mcninch-optimal}*{Prop/Defn 21(4)}, and (4) is
  \cite{mcninch-optimal}*{Prop/Defn 21(5)}.
\end{proof}

Let $\Psi$ be a cocharacter associated with $X$ as in the previous
Proposition. Then the parabolic subgroup $P(X) = P(\Psi)$ of (4) is
known as the the instability parabolic of $X$.

Let $X \in \glie$, and let $[X] \in \Proj(\glie)(K)$ be the $K$-point
which ``is'' the line determined by $X$ in the corresponding projective
space.
\begin{prop}
  \label{prop:N-cochars}
  Write $N_G(X) = \Stab_G([X])$.
  \begin{enumerate}
  \item $N_G(X)$ is a smooth subgroup of $G$.
  \item For each maximal torus $T$ of $N_G(X)$, there is a unique cocharacter
    $\lambda \in X_*(T)$ associated to $X$.
  \end{enumerate}
\end{prop}

\begin{proof}
  Recall that $N_G(X)$ is the \emph{scheme-theoretic} stabilizer of
  the point $[X] \in \Proj(\glie)(K)$.  (1) follows from
  \cite{mcninch-rat}*{Lemma 23}, and in view of Proposition/Definition
  \ref{prop:assoc-cochar}(1), assertion (2) follows from
  \cite{mcninch-rat}*{Lemma 25}.
\end{proof}

\subsection{Notation for $\SL_2$}
\label{sub:sl2-notation}

We fix here some convenient notation for $\SL_2$. We first choose the
``standard'' basis for the Lie algebra $\lie{sl}_2$:
\begin{equation*}
  E = 
  \begin{pmatrix}
    0 & 1 \\
    0 & 0
  \end{pmatrix},
  \quad  
  F =  \begin{pmatrix}
    0 & 0 \\
    1 & 0
  \end{pmatrix},
    \quad \text{and} \quad
    [E,F] = 
  \begin{pmatrix}
    1 & 0 \\
    0 & -1
  \end{pmatrix}.
\end{equation*}  
Now consider the homomorphisms $e,f:\Ga \to \SL_2$ given for each
commutative $K$-algebra $\Lambda$ and each $t \in \Ga(\Lambda) =
\Lambda$ by the rules
\begin{equation*}
  e(t) = 
  \begin{pmatrix}
    1 & t \\
    0 & 1
  \end{pmatrix} 
  \quad \text{and} \quad
  f(t) = 
  \begin{pmatrix}
    1 & 0 \\
    t & 1
  \end{pmatrix}.
\end{equation*}
Finally, write $\T$ for the ``diagonal'' maximal torus of $\SL_2$; we fix the cocharacter
\begin{equation*}
  (t \mapsto     \begin{pmatrix}
  t & 0\\
  0 & t^{-1}
    \end{pmatrix}):\Gm \to \T
\end{equation*}
and use this cocharacter to identify $\T$ with $\Gm$.

\subsection{Optimal homomorphisms}
\label{sub:optimal}
We will use without mention the notation of \S \ref{sub:sl2-notation}.
Let $G$ be a reductive group. We say that a homomorphism $\phi:\SL_2
\to G$ is optimal for $X = d\phi(E)$, or simply that $\phi$ is an
optimal $\SL_2$-homomorphism, if $\lambda = \phi_{\mid \T}$ is a
cocharacter associated with $X$.

\begin{theorem}[\cite{mcninch-optimal}]
  \label{theorem:optimal}
  Suppose that $G$ is strongly standard.  Let $X \in \glie$ satisfy
  $X^{[p]}=0$, and let $\lambda \in X_*(G)$ be associated with $X$.
  There is a unique homomorphism
  \begin{equation*}
    \phi:\SL_2 \to G
  \end{equation*}
  such that $d\phi(E) = X$ and $\phi_{\mid \T} = \lambda$.
  Moreover, the image of $\phi$ is geometrically \gcr.
\end{theorem}

\begin{proof}
  In view of Proposition \ref{prop:strongly-standard-properties}(3),
  this follows from \cite{mcninch-optimal}*{Theorem 47 and Prop. 52}.
\end{proof}

\begin{rem}
  \label{rem:good-vs-optimal}
  Seitz has introduced a notion of ``good $A_1$-subgroup'' of a
  quasisimple group in \cite{seitz}; in \cite{liebeck-seitz}, these
  subgroups are called ``restricted''. Refer to
  \cite{mcninch-optimal}*{\S 8.5} to see that a subgroup of type $A_1$
  of a quasisimple group $G$ is restricted if and only if it is the
  image of an optimal homomorphism $\SL_2 \to G$. It is not hard to
  see that the image of an optimal homomorphism is restricted; cf.
  \cite{mcninch-optimal}*{Prop. 30}. The proof that a restricted
  $A_1$-subgroup is the image of an optimal homomorphism is more involved.
\end{rem}

 \begin{rem}
   \label{rem:SL(V)-optimal}
   If $V$ is a finite dimensional vector space, a homomorphism
   $\phi:\SL_2 \to \SL(V)$ is optimal if and only if $V$ is a
   restricted semisimple $\SL_2$-module. Indeed, if $V$ is restricted
   and semisimple, one sees at once that $\phi_{\mid \T}$ is
   associated with $d\phi(E)$ so that $\phi$ is indeed optimal. On the
   other hand, if $\lambda = \phi_{\mid \T}$ is associated to $X =
   d\phi(E)$, then the the character of the $\SL_2$-module $V$ is
   determined by the cocharacter $\lambda$; it follows that the
   composition factors of $V$ as $\SL_2$-module are restricted.  If $0
   \le n < p$, write $L(n)$ for the restricted simple $\SL_2$-module
   of highest weight $n$ \cite{JRAG}*{\S II.2}. The linkage principle
   \cite{JRAG}*{Corollary II.6.17} implies that
   $\operatorname{Ext}^1_{\SL_2}(L(n),L(m)) =0$ whenever $0 \le n,m <
   p$. Thus, $V$ is semisimple as well.
 \end{rem}

\begin{prop}
  \label{prop:optimal-centralizer}
  Let $S$ be the image of the optimal $\SL_2$-homomorphism $\phi$, and
  let $\lambda = \phi_{\mid \T} \in X_*(G)$. Write 
  $X = d\phi(E)$ and $Y = d\phi(F)$. Then:
  \begin{enumerate}
  \item $C_G(\image d\phi) = C_G(S)$.
  \item The scheme-theoretic intersection $C_G(S) = C_G(X) \cap C_G(\lambda)$ is a smooth subgroup of $G$.
  \end{enumerate}
\end{prop}

\begin{proof}
  For (1), recall that if $\e:\Ga \to G$ is given by $\e = \phi \circ
  e$, then $C_G(X) = C_G(\image \e)$ by \cite{mcninch-optimal}*{Prop.
    35}.  Similarly, $C_G(Y) = C_G(\image \phi \circ f)$. Since
  $\image d\phi$ is spanned by $X$ and $Y$ and since $S$
  is generated as a group scheme by the image of $\phi \circ e$ and
  the image of $\phi \circ f$, we have
  \begin{equation*}
    C_G(\image d\phi) = C_G(X) \cap C_G(Y) = C_G(\image \phi \circ e) 
    \cap C_G(\image \phi \circ f) = C_G(S).
  \end{equation*}

  For (2), the inclusion $C_G(S) \subset C_G(X) \cap C_G(\lambda)$ is
  clear. To prove the other inclusion, let $\Lambda$ be a commutative
  $K$-algebra, and let $g \in G(\Lambda)$ be such that $\Ad(g)X = X$
  and $\Int(g) \circ \lambda = \lambda$. By (1), it is enough to show
  that $g$ centralizes $Y$. Since $Y \in \glie(\lambda;-2)(K)$ and
  since $\Int(g) \circ \lambda = \lambda$, we have $\Ad(g)Y \in
  \glie(\lambda;-2)(\Lambda)$. Notice that
  \begin{equation*}
    [X,\Ad(g)Y] = [\Ad(g^{-1})X,Y] = [X,Y].
  \end{equation*}
  Thus $[X,Y-\Ad(g)Y] = 0$ so that $Y-\Ad(g)Y \in
  \clie_\glie(X)(\lambda;-2)(\Lambda)$.  Since $\clie_\glie(X) \subset
  \Lie(P(\lambda))$ by Proposition \ref{prop:assoc-cochar}, we have
  $\clie_\glie(X)(\lambda;-2) = 0$ so that $Y=\Ad(g)Y$ as
  required. Now, $C_G(X)$ is smooth by Proposition
  \ref{prop:strongly-standard-properties}, hence $C_G(X) \cap C_G(\lambda)$ is
  smooth by \cite{SGA3}*{II Exp. XI, Cor 5.3}.

\end{proof}

\begin{rem}
  \label{rem:sl2-image}
  In the notation of the previous Proposition, we have $\image d\phi =
  \Lie S$ whenever $p>2$, since the adjoint representation of $\SL_2$
  is irreducible for $p>2$.
\end{rem}

\begin{prop} Let $G$ be a strongly standard reductive group.
  \label{prop:optimal-in-Levi-and-under-isogeny}
  \begin{enumerate}
  \item   Let $L \subset G$ be a Levi subgroup, and assume that $\phi:\SL_2 \to L$
    is a homomorphism. Then $\phi$ is an optimal homomorphism in $G$ if
    and only if it is an optimal homomorphism in $L$.
  \item Let $\pi:G_1 \to G$ be a central isogeny, let $f:\SL_2 \to G$
    be a homomorphism, and suppose that $\widetilde{f}:\SL_2 \to G_1$
    satisfies $\pi \circ \widetilde{f} = f$. Then $f$ is optimal if
    and only if $\widetilde{f}$ is optimal.
  \end{enumerate}
\end{prop}

\begin{proof}
  We first prove (1). In view of Theorem \ref{theorem:optimal}, it
  suffices to prove the following: Let $X \in \Lie(L)$ be nilpotent
  and let $\lambda \in X_*(L)$ be a cocharacter with $X \in
  \Lie(L)(\lambda;2)$. Then $\lambda$ is associated to $X$ in $L$ if
  and only if $\lambda$ is associated to $X$ in $G$.

  Note that $\lambda \in X_*(N_L(X))$, so the image of $\lambda$
  normalizes $C_L(X)$. In particular, we may choose a maximal torus
  $S_0$ of $C_L(X)$ centralized by the image of $\lambda$, and we
  may choose a maximal torus $S$ of $C_G(X)$ with $S_0 \subset S$.
  Notice that $S_0$ -- and hence also $S$ -- contains the center of
  $L$.  Since $L$ is the centralizer in $G$ of the connected center of $L$,
  we have $S \subset L$ so that $S=S_0$.  Moreover, since $C_G(S)
  \subset L$, it is clear that $C_G(S) = C_L(S) = M$.  Moreover, it is
  clear that $\lambda \in X_*(M)$, and the Proposition follows since
  the condition that $\lambda$ be associated to $X$ is just that
  $\lambda \in X_*((M,M))$; this condition is the same for $L$ and for
  $G$.

  We now prove (2). Note first that it suffices to prove (2) in case
  $K$ is algebraically closed. Let $X = df(E)$ and
  $\widetilde{X} = d\widetilde{f}(E)$.  Then $\pi$ induces a
  surjective morphism $N_{G_1}(\widetilde{X})_\red \to N_G(X)$.
  \footnote{Note that $G_1$ need not be strongly standard.} We may
  thus choose a maximal torus $\widetilde{T}$ of 
  $C_{G_1}(\widetilde{X})_\red$ centralized by $\image
  \widetilde{f}_{\mid \T}$ and a maximal torus $T$ of $C_G(X)$
  centralized by $\image f_{\mid T}$ such that $\pi(\widetilde{T}) =
  T$.

  Now, $\widetilde{X}$ is distinguished in the Levi subgroup $L_1 =
  C_{G_1}(\widetilde{T})$ and $X$ is distinguished in the Levi
  subgroup $L = C_G(T)$. Since the maximal tori in
  $C_{G_1}(\widetilde{X})$ are all conjugate, one sees that
  $\widetilde{f}_{\mid \T}$ is associated with $\widetilde{X}$ if and
  only if $\image \widetilde{f}_{\mid \T}$ lies in the derived group of $L_1$;
  since the maximal tori in $C_G(X)$ are all conjugate, $f_{\mid \T}$ is
  associated with $X$ if and only if $\image f_{\T}$ lies in the derived
  group of $L$. Since $\pi$ induces a central isogeny $L_1 \to L$, it
  follows that $\widetilde{f}_{\mid T}$ is associated with
  $\widetilde{X}$ if and only if $f_{\mid T}$ is associated with $X$;
   (2) is an immediate consequence.
\end{proof}

\subsection{Frobenius endomorphisms}
\label{sub:frob}
Let $H$ be a connected, split, quasi-simple algebraic group; recall
that $H$ arises by base change from a corresponding group scheme
$H_{/\F_p}$ over the prime field $\F_p$.  There is a Frobenius
endomorphism $F:H \to H$ which arises by base change from the
corresponding Frobenius endomorphism of $H_{/\F_p}$.
\begin{prop}
  \label{prop:frobenius} Let $G$ be an algebraic group, and let
  $\phi:H \to G$ be a homomorphism. The following are equivalent:
  \begin{enumerate}
  \item $d\phi = 0$
  \item there is a unique integer $t \ge 1$ and a unique homomorphism
    $\psi:H \to G$ such that $d\psi \ne 0$ and $\phi = \psi
    \circ F^t$.
  \end{enumerate}
\end{prop}

\begin{proof}
  (1) $\implies$ (2) is a consequence of \cite{mcninch-optimal}*{Cor.
    20}.
  (2) $\implies$ (1) is straightforward.
\end{proof}




Of course, the above result hold in particular when $H$ is the group $\SL_2$.

\section{The tangent map of a $G$-completely reducible $\SL_2$-homomorphism}
\label{sec:lie-optimal}

\subsection{The set-up}
Now fix a homomorphism $\phi:\SL_2 \to G$ whose image is geometrically
\gcr.  Assume that $d\phi \ne 0$, and write
\begin{equation*}
  X = d\phi(E), \quad  Y = d\phi(F), \quad H = d\phi([E,F]) \in \glie,
\end{equation*}
Also put $\slie = \image d\phi$, and write $\lambda = \phi_{\mid \T}$.
Consider the smooth subgroups $N_G(X), N_G(Y) \subset G$ which are the stabilizers
of the points $[X], [Y] \in \Proj(\glie)(K)$. Then $\lambda$
is evidently a cocharacter of $N_G(X) \cap N_G(Y)$.

Consider the group schemes $C(X,Y) = (C_G(X) \cap C_G(Y))$ and $N(X,Y) = (N_G(X)
\cap N_G(Y))$.  We observe the following:
\begin{lem}
   $C(X,Y)$ is a normal subgroup scheme of $N(X,Y)$.
\end{lem}
In particular, the image of $\lambda$ normalizes $C(X,Y)$.

\subsection{Working geometrically}
\label{sub:geom}
Fix an algebraically closed field $k$ containing $K$ and consider
$G_{/k}$, $\slie_{/k} = \slie \tensor_K k$ etc.  In this section, we
are forced to consider the reduced subgroups corresponding to various
subgroup schemes; recall the results of \ref{sub:reduced}. Thus,
\emph{for the remainder of \S\ref{sub:geom}, we replace $K$ by $k$ and
  so suppose that $K$ is algebraically closed.}

According to \S \ref{sub:reduced}, the image of $\lambda$ normalizes
$C(X,Y)$ and hence also $C(X,Y)_\red$. Thus, we may choose a maximal
torus $T \subset C(X,Y)_\red$ centralized by the image of $\lambda$.

Consider now the Levi subgroup $M = C_G(T)$ of $G$; $M$ is a strongly
standard reductive group by Proposition
\ref{prop:strongly-standard-properties}(1).  Since $X$ and $Y$ are
centralized by $T$, and since $\slie$ is generated as a Lie algebra by
$X$ and $Y$, we have $\slie \subset \Lie(M)$. Of course, the image of
the homomorphism $\phi$ need not lie in $M$.

\begin{lem}
  \label{lem:slie-not-in-P}
  $\slie$ is not contained in $\Lie(P)$ for any proper parabolic
  subgroup $P \subset M$.
\end{lem}

\begin{proof}
  Any torus $T_1 \subset M$ centralizing $\slie$ of course centralizes $X$
  and $Y$; thus $T$ lies in a maximal torus of $C(X,Y)_\red$. Since
  $T$ is central in $M$, $T_1$ centralizes $T$. Since $T$ is a maximal
  torus of $C(X,Y)_\red$, we find that $T_1 \subset  T$ hence $T_1$ is
  central in $M$.  Since $\slie$ is the Lie algebra of a $G$-cr
  subgroup of $G$, the Lie algebra $\slie$ is itself $G$-cr by Theorem
  \ref{theorem:lie-cr}. Hence $\slie$ is also $M$-cr by Proposition
  \ref{prop:gcr-levi} . If $\slie$ is contained in $\Lie(P)$ for a
  parabolic subgroup $P\subset M$, then $\slie$ is contained in $\Lie(L)$ for
  some Levi subgroup $L$ of $P$. But then any central torus of $L$ is
  central in $M$, so that $P=M$.
\end{proof}

\begin{prop}
  \label{prop:lambda0-lambda-centralizes-slie}
  Let $T_1$ be a maximal torus of $N_M(X)$ with $\lambda \in X_*(T_1)$,
  and let $\lambda_0 \in X_*(T_1)$ be the unique cocharacter of $T_1$
  associated to $X$ [see Proposition \ref{prop:N-cochars}(2)]. Let
  $\phi_0:\SL_2 \to M$ be the optimal homomorphism determined by $X$
  and $\lambda_0$ [Theorem \ref{theorem:optimal}].  Write $\mu =
  \lambda_0 - \lambda$ for the cocharacter
  \begin{equation*}
    t \mapsto    \lambda_0(t)\cdot\lambda(t^{-1})   
  \end{equation*}
  of $T_1$. Then the image of $\mu$ is central in $M$.
\end{prop}

\begin{proof}
  We have $H \in \mlie(\lambda_0;0) \cap \mlie(\lambda;0)$ by the
  choice of $T_1$; thus $H \in \mlie(\mu;0)$. We have also
  $X \in \mlie(\lambda_0;2)  \cap \mlie(\lambda;2)$ so that $X
  \in \mlie(\mu;0)$ as well.

  Write $Y = \sum_{j \in \Z} Y^j$ with $Y^j \in \mlie(\lambda_0;j)$.
  Since $[X,Y] = H \in \mlie(\lambda_0;0)$, we have $Y-Y^{-2} \in
  \clie_{\mlie}(X)$, so that
  \begin{equation*}
    Y = Y^{-2} + \sum_{j \ge 0} Y^j.
  \end{equation*}
  Since the images of $\lambda$ and $\lambda_0$ commute, and since $Y
  \in \mlie(\lambda;-2)$, we have $Y^j \in \mlie(\lambda;-2)$ for all
  $j$.  Thus, $Y^j \in \mlie(\lambda_0 - \lambda;j+2)=\mlie(\mu;j+2)$
  for all $j$, hence $Y \in \sum_{\ell \ge 0} \mlie(\mu;\ell) = \Lie
  P_M(\mu)$.

  Since $X,Y,H \in \Lie P_M(\mu)$, we have proved that $\slie=\image
  d\phi$ lies in $\Lie P_M(\mu)$.  Thus by Lemma
  \ref{lem:slie-not-in-P}, we have $P_M(\mu) = M$; we conclude that
  the image of $\mu$ is central in $M$.  
\end{proof}

\begin{prop}
    \label{prop:d-phi-equals-d-phi-0}
    Let $T_1$ be a maximal torus of $N_M(X)$, and write $\phi_0:\SL_2
    \to M$ for the optimal homomorphism determined by the cocharacter
    $\lambda_0 \in X_*(T_1)$ associated with $X$ as in Proposition
    \ref{prop:lambda0-lambda-centralizes-slie}.  Then $d\phi =
    d\phi_0$.
\end{prop}

\begin{proof}
  Recall that $T$ is a fixed maximal torus of $C(X,Y)_\red$, and $M=
  C_G(T)$.  Using \eqref{eq:cover-G}, one finds a (possibly
  inseparable) central isogeny
  \begin{equation*}
    \pi:G_\sconn \to G,
  \end{equation*}
  where the derived group of $G_\sconn$ is simply connected.  There is
  a torus $\widetilde{T} \subset G_\sconn$ with $\pi(\widetilde{T}) = T$;
  then the Levi subgroup $M_\sconn = C_{G_\sconn}(\widetilde{T})$ has
  simply connected derived group, and $\pi$ restricts to a central
  isogeny $\pi:M_\sconn \to M$.

  Since $\SL_2$ is simply connected, there are homomorphisms
  \begin{equation*}
    \widetilde{\phi}:\SL_2 \to G_\sconn \quad \text{and} 
    \quad \widetilde{\phi}_0:\SL_2
    \to M_\sconn   
  \end{equation*}
  such that $\phi = \pi \circ \widetilde{\phi}$ and
  $\phi_0 = \pi \circ \widetilde{\phi}_0$. We write
  $\widetilde{\lambda}$ and $\widetilde{\lambda}_0$ for the
  cocharacters obtained by restricting these homomorphisms to the
  maximal torus $\T$ of $\SL_2$.  Moreover, write
  \begin{equation*}
    \mu = \lambda_0 - \lambda_1 \quad \text{and} \quad
    \widetilde{\mu} =   \widetilde{\lambda}_0 -   \widetilde{\lambda}_1.
  \end{equation*}
  Since $\lambda,\lambda_0 \in X_*(T_1)$, we know that
  $\widetilde{\lambda},\widetilde{\lambda}_0$ are cocharacters of
  $M_\sconn$.  Since $\pi \circ \widetilde{\mu} = \mu$, it follows
  from Proposition \ref{prop:lambda0-lambda-centralizes-slie} that the
  image of $\widetilde{\mu}$ is central in $M_\sconn$.

  Since $T$ centralizes $\image \phi_0$ and $\image d\phi$, it is
  clear that the images of $d\widetilde{\phi}_0$ and
  $d\widetilde{\phi}$ lie in $\Lie(M_\sconn)$.  Write $H_0 =
  d\phi_0([E,F])$. We claim that the Proposition will follow if we
  show that $H=H_0$.  Indeed, by Proposition
  \ref{prop:lambda0-lambda-centralizes-slie}, the image of $\mu =
  \lambda_0 - \lambda$ centralizes $Y$; thus, we have $Y \in
  \mlie(\lambda_0;-2)$ and $Y-Y_0 \in \mlie(\lambda_0;-2)$. If
  $H=H_0$, then $[X,Y-Y_0] = H - H_0 = 0$ so that $Y-Y_0 \in
  \clie_\mlie(X)$. Since $\lambda_0$ is associated to $X$, we have
  $\clie_\mlie(X) \subset \Lie P(\lambda_0) = \sum_{i \ge 0}
  \mlie(\lambda_0;i)$ by Proposition/Definition
  \ref{prop:assoc-cochar}.  We may thus conclude that $Y=Y_0$ so that
  $d\phi$ and $d\phi_0$ are indeed equal.

  It remains now to show that $H=H_0$. Let $\widetilde{H} =
  d\widetilde\phi([E,F])$ and $\widetilde{H}_0 =
  d\widetilde\phi_0([E,F])$.  It is clearly enough to show that
  $\widetilde{H} = \widetilde{H}_0$. 

  Now, since the derived group $M_\sconn'$ of $M_\sconn$ is simply
  connected, one knows that $M_\sconn$ is a direct product
  \begin{equation*}
    M_\sconn = Z^o(M_\sconn) \times M_\sconn'
  \end{equation*}
  where $Z^o(M_\sconn)$ is the connected component of the center of
  $M_\sconn$ (it is a torus).
  Thus also
  \begin{equation*}
    \Lie(M_\sconn) = \Lie(Z^o(M_\sconn)) \oplus \Lie(M_\sconn').
  \end{equation*}
  
  The derived subalgebra $[\Lie(M_\sconn),\Lie(M_\sconn)]$ is
  contained in $\Lie(M_\sconn')$; since $\sllie_2 =
  [\sllie_2,\sllie_2]$, we have
  \begin{equation*}
    \image d\widetilde{\phi} \subset \Lie(M_\sconn') \quad
    \text{and}
    \quad \image d\widetilde{\phi}_0 \subset \Lie(M_\sconn').
  \end{equation*}
  In particular,
  \begin{equation}
    \label{eq:H-H0-in-Lie(M')}
    \widetilde{H}_0 - \widetilde{H} \in \Lie(M_\sconn').
  \end{equation}
  On the other hand, $\image \widetilde{\mu}$ lies in $Z^o(M_\sconn)$,
  and so 
  \begin{equation}
    \label{eq:H-H0-in-Lie(Z)}
    \widetilde{H}_0 - \widetilde{H} \in \image d\widetilde{\mu}
  \subset \Lie(Z^o(M_\sconn)).
  \end{equation}
  Since $\Lie(Z^o(M_\sconn)) \cap \Lie(M_\sconn') = 0$, we deduce that
  $\widetilde{H}_0 = \widetilde{H}$
  by applying \eqref{eq:H-H0-in-Lie(M')} and
  \eqref{eq:H-H0-in-Lie(Z)}. This completes the proof.
\end{proof}

\subsection{The tangent map over any field}
\label{sub:ground-field}

We now suppose that $K$ is an arbitrary field of characteristic
$p>0$.  As in the previous section, we fix a homomorphism $\phi:\SL_2
\to G$ whose image is geometrically \gcr, and we assume that $d\phi
\ne 0$.

We also fix an algebraically closed field $k$ containing $K$.
\begin{cor}
  \label{cor:smooth-intersection}
  The $K$-subgroup $C(X,Y) = C_G(X) \cap C_G(Y)$ is smooth.
\end{cor}

\begin{proof}
  Working over the algebraically closed field $k \supset K$, let
  $\phi_0:\SL_{2/k} \to G_{/k}$ be any optimal $k$-homomorphism as in
  Proposition \ref{prop:d-phi-equals-d-phi-0}. Write $S_0 \subset G_{/k}$
  for the image of $\phi_0$, and recall that $\slie_{/k} = \image
  d\phi_{/k} = \image d\phi_0$.  Then we know that
  \begin{equation*}
    C_{G_{/k}}(S_0) = C_{G_{/k}}(\slie_{/k}) = C_{G_{/k}}(X) \cap C_{G_{/k}}(Y)
  \end{equation*}
  by Proposition \ref{prop:optimal-centralizer}, hence $C_{G_{/k}}(X)
  \cap C_{G_{/k}}(Y) = (C_G(X) \cap C_G(Y))_{/k}$ is smooth. But then
  $C_G(X) \cap C_G(Y)$ is smooth, since that is so after extension of
  the ground field.
\end{proof}

\begin{cor}
  \label{cor:gcr-tangent}
  There is a cocharacter $\lambda_0$ of $G$ associated to $X$ such
  that if $\phi_0:\SL_2 \to G$ is the optimal homomorphism determined
  by $X$ and $\lambda_0$, then $d\phi = d\phi_0$. Moreover, $\phi_0$
  is uniquely determined by $\phi$.
\end{cor}

\begin{proof}
  Since $C_G(X) \cap C_G(Y)$ is smooth by the previous corollary, we
  can find a maximal torus $T$ of $C_G(X) \cap C_G(Y)$ centralized by the
  image of the torus $\lambda = \phi_{\mid \T}$. Then the Levi
  subgroup $M=C_G(T)$ is a strongly standard reductive $K$-subgroup.
  As in Proposition \ref{prop:d-phi-equals-d-phi-0}(1), we may find
  maximal tori [now over $K$] of $C_M(X)$ centralized by the image of
  $\lambda$; Proposition \ref{prop:d-phi-equals-d-phi-0} then gives
  the first assertion of the Corollary. According to Theorem
  \ref{theorem:optimal}, an optimal homomorphism is uniquely
  determined by its tangent mapping; the uniqueness assertion follows
  at once.
\end{proof}

\section{Proof of the main theorem}
\label{sec:main-theorem}

\subsection{A general setting}
\label{sub:setting}
Let $H$ be a connected and simple algebraic group.  For each strongly
standard reductive group $G$, suppose that one is given a set $\C_G$
of homomorphisms $H \to G$ with the properties to be enumerated below.

Let $G$ be strongly standard and let $f_0 \in \C_G$ be arbitrary; 
write $S_0$ for the image of $f_0$.  We assume the following hold for
each $f_0$:
\begin{enumerate}
\item[(C1)] $S_0$ is geometrically \gcr.
\item[(C2)] $C_G(S_0)$ is a smooth subgroup of $G$, and $C_G(S_0) =
  C_G(\Lie(S_0))$.
\item[(C3)] $\Lie(S_0) = \image df_0$.
\end{enumerate}
We also suppose:
\begin{enumerate}
\item[(C4)] Given any homomorphism $f:H \to G$ for which $df \ne 0$
  and for which $\image f$ is geometrically \gcr, there is a unique
  $f_0 \in \C_G$ such that $df = df_0$.
\item[(C5)] If $f:H \to G$ is a homomorphism and if $L \subset G$ is a Levi subgroup
   with $\image f \subset L$, then $f \in \C_G$ if and only if $f
  \in \C_L$.
\end{enumerate}

The following Lemma gives a useful application of (C1) and (C2).
\begin{lem}
  \label{lem:parabolic-inclusion}
  Let $G$ be a reductive group and let $S \subset G$ be a subgroup with the
  property $C_G(S) = C_G(\Lie(S))$. Suppose that $S$ is geometrically
  \gcr.  If $K \subset k$ is any field extension and $P \subset G_{/k}$ is a
  $k$-parabolic subgroup, then $S_{/k} \subset P$ if and only if
  $\Lie(S)_{/k} \subset \Lie(P)$,
\end{lem}

\begin{proof}
  Since the Lemma follows once it is proved for algebraically closed
  extensions $k$, it suffices to suppose that $K$ itself is
  algebraically closed and to prove the conclusion of the Lemma for a
  ($K$-)parabolic subgroup $P \subset G$. First notice that if $S \subset P$,
  then clearly $\Lie(S) \subset \Lie(P)$.

  Now suppose that $\slie = \Lie(S) \subset \Lie(P)$.  Since $S$ is
  \gcr, Theorem \ref{theorem:lie-cr} shows that $\slie$ is also \gcr.
  Thus, we may find a Levi subgroup $L \subset P$ with $\slie \subset
  \Lie(L)$.  Then $L = C_G(T)$ where $T = Z(L)$, and we see that
  \begin{equation*}
    T \subset  C_G(\slie) = C_G(S);
  \end{equation*}
  thus $T$ centralizes $S$, so that $S \subset C_G(T) = L \subset P$, as required.
\end{proof}

We now observe:
\begin{prop}
  \label{prop:sl2-conditions-hold}
  Let $p>2$, let $H = \SL_2$, and for each strongly standard reductive
  group $G$, let $\C_G$ be the set of optimal homomorphisms $\SL_2 \to
  G$. Then conditions (C1) -- (C5) of \S \ref{sub:setting} hold for
  the sets $\C_G$.
\end{prop}

\begin{proof}
  (C1) follows from Theorem \ref{theorem:optimal}, (C2) is Proposition
  \ref{prop:optimal-centralizer}, (C4) is Corollary
  \ref{cor:gcr-tangent}, and (C5) is Proposition
  \ref{prop:optimal-in-Levi-and-under-isogeny}(1).
  
  Since $p>2$, the adjoint representation of $\SL_2$ is irreducible;
  since any optimal homomorphism $f:\SL_2 \to G$ has $df(E) \ne 0$,
  the map $df$ must be injective and so (C3) is immediate.
\end{proof}

\subsection{Some results about twisted-product homomorphisms}
\label{sub:twisted-product-results}

Let $H$ be a reductive group and let $\C_G$ be a collection of
homomorphisms $H \to G$ for each strongly standard group $G$ which
satisfies (C1)--(C5) of \S \ref{sub:setting}.

We are going to prove several technical results about twisted product
homomorphisms; to avoid repetition in the statements, we fix the
following notation:

Let $\vec h = (h_0,h_1,\dots,h_r)$ with $h_i \in \C_G$ be commuting
homomorphisms [as in the introduction], and let $\vec n = (n_0 < n_1 <
\cdots < n_r)$ be non-negative integers; the data $(\vec h, \vec n)$
determines a twisted-product homomorphism $\Phi = \Phi_{\vec h,\vec
  n}:H \to G$ given for each commutative $K$-algebra $\Lambda$ and
each $g \in H(\Lambda)$ by the rule
\begin{equation}
  \label{eq:twisted-hom}
  g \mapsto h_0(F^{n_0}g) \cdot h_1(F^{n_1}g) \cdots h_r(F^{n_r}g).
\end{equation}

Several of the results proved in this section hold only assuming a
subset of the conditions (C1)--(C5); for simplicity of exposition, we
assume all five conditions hold -- we don't bother to identify the
 subset.

\begin{lem}
  \label{lem:untwisted-product}
  Let $G$ be strongly standard, let $\vec h$, $\vec n$, and
  $\Phi=\Phi_{\vec h,\vec n}$ be as in the beginning of \S
  \ref{sub:twisted-product-results}. Then $d\Phi = 0$ if and only if
  $n_0>0$. If $d\Phi = 0$ let $\Psi$ be the twisted-product
  homomorphism determined by $(\vec h, (0,n_1 - n_0,\dots,n_r-n_0))$.
  Then $\Phi = \Psi \circ F^{n_0}$ and $d\Psi \ne 0$.  Moreover,
  $\image \Phi = \image \Psi$.
\end{lem}

\begin{proof}
  Straightforward and left to the reader.
\end{proof}

\begin{prop}
  \label{prop:twisted-product-gcr}
  Let $G$ be strongly standard, let $(f_1,f_2)$ be commuting
  homomorphisms $H \to G$ with $f_1 \in \C_G$. Let $(n_1 < n_2)$ be
  non-negative integers, and let $f$ be the twisted-product
  homomorphism determined by $(f_1,f_2)$ and $(n_1,n_2)$.  Write $S_i$
  for the image of $f_i$, $i=1,2$, and write $S$ for the image of $f$.
  Then:
  \begin{enumerate}
  \item for each field extension $K \subset k$ and each parabolic
    subgroup $P \subset G_{/k}$, we have $S_{/k} \subset P$ if and only if
    $S_{1/k}\subset P$ and $S_{2/k} \subset P$.
  \item $S_1 \cdot S_2$ is geometrically \gcr\ if and only if $S$ is
    geometrically \gcr.
  \end{enumerate}
\end{prop}

\begin{proof}
  Note that (1) will follow once it is proved for algebraically closed
  extension fields $k$ of $K$.  Thus, we suppose that $k=K$ is
  algebraically closed, and prove the conclusion of (1) for parabolic
  subgroups $P \subset G$.

  If the parabolic subgroup $P \subset G$ contains $S_1$ and $S_2$, it is
  clear by the definition of a twisted-product homomorphism that $P$
  contains $S$. Suppose now that $P$ contains $S$; we show it contains
  also $S_1$ and $S_2$.

  Applying Lemma \ref{lem:untwisted-product}, one knows that if $g:H
  \to G$ is the twisted-product homomorphism determined by $(f_1,f_2)$
  and $(0,n_2 - n_1)$, then $\image g = S$ as well. We may thus
  suppose that $n_1 = 0$, so that $df \ne 0$.

  It is clear that $df = df_1$. Since $\image df_1 = \Lie(S_1)$, it
  follows that $\Lie(S) = \Lie(S_1)$. Since $S \subset P$, we have
  $\Lie(S_1) = \Lie(S) \subset \Lie(P)$; since (C1) and (C2) hold, we may
  apply Lemma \ref{lem:parabolic-inclusion} and conclude that $S_1 \subset
  P$.  Since $f_2$ is given by the rule
  \begin{equation*}
    g \mapsto f_1(g)^{-1}f(g)
  \end{equation*}
  it is then clear that $S_2 \subset P$ as well. This proves (1).

  Since (2) is a geometric statement, we may again suppose that $K$ is
  an algebraically closed field.  Write $X$ for the building of $G$;
  cf.  \cite{serre-sem-bourb}*{\S 2 and \S 3.1}.  Then $X$ is a
  simplicial complex whose simplices are in bijection with the
  parabolic subgroups of $G$. We have shown the equality of
  fixed-point sets: $X^S = (X^{S_2})^{S_1} = X^{S_1 \cdot
    S_2}$.

  According to \cite{serre-sem-bourb}*{Th\'eor\`eme 2.1}, the group
  $S$ is \gcr\ if and only if $X^S= X^{S_1S_2}$ is contractible if and only if $S_1
  \cdot S_2$ is \gcr. This proves (2).
\end{proof}

\begin{cor}
  \label{cor:twisted-product-in-P}
  Let $G$ be strongly standard, and let $\vec h$, $\vec n$, $\Phi =
  \Phi_{\vec h,\vec n}$ be as in the beginning of \S
  \ref{sub:twisted-product-results}.  Write $S$ for the image of $h$
  and $S_i$ for the image of $h_i$.  If $P \subset G_{/k}$ is a
  $k$-parabolic subgroup for an  extension field
  $K \subset k$, then $S_{/k} \subset P$ if and only if $S_{i/k}
  \subset P$ for $i=1,\dots,r$.
\end{cor}

\begin{proof}
  It is enough to give the proof assuming that $K=k$ is algebraically
  closed.  If $S_i \subset P$ for each $i$, it is clear by construction that
  $S \subset P$.  

  Now suppose that $S \subset P$. To prove that each $S_i \subset P$, we proceed by
  induction on $r$. If $r=1$, the result is immediate.  Suppose that
  $r>1$, and let $\Psi:H \to G$ be the twisted-product homomorphism
  determined by $(h_2,\dots,h_r)$ and $(n_2 - 1,\dots,n_r - 1)$.  Then
  $f$ may be regarded as the twisted-product homomorphism determined
  by $(h_1,\Psi)$ and $(0,1)$. Thus we may
  apply Proposition \ref{prop:twisted-product-gcr}(1) to see that $S_1
  \subset P$ and $\image \Psi \subset P$. Now apply the induction hypothesis to
  $\Psi$ to learn that $S_i \subset P$ for $2 \le i \le r$. This completes
  the proof.
\end{proof}

\begin{prop}
  \label{prop:twisted-centralizer}
  Let $G$ be strongly standard, and let $\vec h$, $\vec n$, $\Phi =
  \Phi_{\vec h,\vec n}$ be as in the beginning of \S
  \ref{sub:twisted-product-results}.  If $S$ denotes the image of $\Phi$
  and $S_i$ the image of $h_i$, then $C_G(S) \subset C_G(S_i)$ for $1
  \le i \le r$.
\end{prop}

\begin{proof}
  If $r=1$, the result is immediate. Suppose $r > 1$, write $\Psi$ for
  the homomorphism determined by $(h_2,\dots,h_r)$ and
  $(n_2,\dots,n_r)$ and write $T = \image \Psi$.  It suffices by
  induction on $r$ to show that $C_G(S) \subset C_G(S_1)$ and $C_G(S)
  \subset C_G(T)$, since then for $2 \le i \le r$ we have
  \begin{equation*}
    C_G(S) \subset C_G(T) \subset C_G(S_i)
  \end{equation*}
  by the induction hypothesis.
  
  Applying Lemma \ref{lem:untwisted-product}, we may assume that $n_1
  = 0$ and $dh \ne 0$ without changing $S$.  Thus $\Lie(S) =
  \Lie(S_1)$.  By (C2), we have
  \begin{equation*}
    C_G(S) \subset C_G(\Lie(S))  = C_G(\Lie(S_1)) = C_G(S_1).
  \end{equation*}

  Finally, it remains to check that $C_G(S) \subset C_G(T)$.  Write
  $\Psi^*,h_1^*,\Phi^*:K[G] \to K[H]$ for the comorphisms of $\Psi$,
  $h_1$, and $\Phi$.  Then by construction, $\Psi^*$ is given by the
  composition
  \begin{equation*}
    K[G] \xrightarrow{\mu} K[G] \tensor_K K[G] \xrightarrow{\iota \tensor id} 
    K[G] \tensor_K K[G] \xrightarrow{h_1^* \tensor \Phi^*} K[H] \tensor_K K[H] 
    \xrightarrow{\Delta} K[H]
  \end{equation*}
  where the map $\mu$ defines the multiplication in $G$, the map
  $\iota$ defines the inversion in $G$, and $\Delta$ is given by
  multiplication in $K[H]$.

  Let $g \in C_G(S)(\Lambda)$ for some commutative $K$-algebra
  $\Lambda$. To show that $g \in C_G(T)(\Lambda)$, it is enough to
  argue that the inner automorphism $\Int(g)$ of $G$ induces the
  identity on the subgroup scheme $T_{/\Lambda}$. Since $T$ is defined
  by the ideal $\ker \Psi^* \lhd K[G]$, it is enough to require that
  $\Psi^*(\Int(g)^*f)= \Psi^*(f)$ for each $f \in \Lambda[G]$.  [Note:
  we write $\Psi^*$ rather than $\Psi^*_{/\Lambda}$ for simplicity.]

  Since $g \in C_G(S)(\Lambda)$ and $g \in C_G(S_1)(\Lambda)$, we know
  for each $f \in \Lambda[G]$ that
  \begin{equation*}
    h_1^*(\Int(g)^*f) = h_1^*(f) \quad \text{and} 
    \quad \Phi^*(\Int(g)^*f) = \Phi^*(f).
  \end{equation*}
  If $f_1 \tensor f_2 \in \Lambda[G] \tensor_\Lambda \Lambda[G]$, then
  \begin{align*}
        (h_1^* \tensor \Phi^*)((\Int(g)^*\tensor \Int(g)^*)(f_1 \tensor f_2)) 
        &=  h_1^*(\Int(g)^*f_1)\tensor \Phi^*(\Int(g)^*f_2)  \\
        &=  h_1^*(f_1)\tensor \Phi^*(f_2)  \\
        &=  (h_1^* \tensor \Phi^*)(f_1 \tensor f_2).
  \end{align*}
 It follows for any $f_1 \in \Lambda[G] \tensor_\Lambda \Lambda[G]$ that
 \begin{equation}
   \label{eq:part-1}
   (h_1^* \tensor \Phi^*)((\Int(g)^* \tensor \Int(g)^*)f_1) = (h_1^* \tensor \Phi^*)f_1
 \end{equation}

  Since $\Int(g)$ is an automorphism of $G$, we have for $f \in
  \Lambda[G]$ that
  \begin{equation}
    \label{eq:part-2}
    ((\iota \tensor id) \circ \mu)(\Int(g)^*f) 
    = (\Int(g)^* \tensor \Int(g)^*) ((\iota \tensor id) \circ \mu)(f).
  \end{equation}

  Combining \eqref{eq:part-1} and \eqref{eq:part-2}, we see that
  $\Psi^*(\Int(g)^*f) = \Psi^*(f)$ for each $f \in \Lambda[G]$, as
  required. This completes the proof.
\end{proof}

\begin{rem}
  With notation as before, one can even show that
  \begin{equation*}
    C_G(S) = \bigcap_{i=1}^r C_G(S_i).
  \end{equation*}
  The inclusion $C_G(S) \subset \bigcap_i C_G(S_i)$ follows from the
  previous Proposition, and the reverse inclusion may be proved by
  showing for each commutative $K$-algebra $\Lambda$ that if $g \in
  C_G(S_i)(\Lambda)$ for each $i$, then $\Phi^*(\Int(g)^*f) = \Phi^*(f)$ for
  each $f \in \Lambda[G]$; the proof is like that used for the
  Proposition.
\end{rem}

\subsection{Finding the twisted factors of a homomorphism with \gcr\ image}
\label{sub:proof-in-general}

Let $H$ be a reductive group and let $\C_G$ be a collection of
homomorphisms $H \to G$ for each strongly standard group $G$ which
satisfies (C1)--(C5) of \S \ref{sub:setting}.
In this section, we are going to give the proof of Theorem \ref{theorem:main}.

We first have the following:
\begin{prop} 
  \label{prop:steinberg-untwist}
  Fix a strongly standard reductive group $G$, and let the
  homomorphism $f:H \to G$ have geometrically \gcr\ image $S$.  Assume
  that $df \ne 0$ and let $f_0 \in \C_G$ be the unique map -- as in
  (C4) -- such that $df = df_0$.  Then:
  \begin{enumerate}
  \item the map $f_1: H \to G$ given by the rule $g \mapsto
    f_0(g^{-1})\cdot f(g)$ is a group homomorphism.
  \item $S_1 = \image f_1 \subset C_G(\image f_0)$.
  \item $df_1 = 0$.
  \item $S_1$ is geometrically \gcr.
  \end{enumerate}
\end{prop}

\begin{proof}
  Write $S_0 = \image f_0$, and write $f_1:H \to G$ for the morphism
  defined by the rule in (1).  Let $\Lambda$ be an arbitrary
  commutative $K$-algebra, let $g \in H(\Lambda)$, and let $X \in
  \Lie(H)(\Lambda)$.   Since $df  = df_0$, we know that
  \begin{equation*}
    \Ad(f(g))df_0(X) = \Ad(f(g))df(X) = df(\Ad(g)X)  = df_0(\Ad(g)X) = 
    \Ad(f_0(g))df_0(X).
  \end{equation*}
  It follows that $\Ad(f_1(g)) = \Ad(f_0(g^{-1})f(g))$ centralizes
  $df_0(X)$ for each $X \in \Lie(H)(\Lambda)$. Since $\image df_0 =
  \Lie(S_0)$ by (C3), it follows that the image of $f_1$ lies in
  $C_G(\Lie(S_0))$.  If now $g,h \in H(\Lambda)$, then we see that
  \begin{equation*}
    f_1(gh)  = f_0(h^{-1})f_0(g^{-1})f(g)f(h) = f_0(h^{-1})f_1(g)f(h) =
    f_1(g)f_0(h^{-1})f(h) = f_1(g)f_1(h).
  \end{equation*}
  Thus $f_1$ is a homomorphism, so that (1) and (2) are proved.
  By construction, the tangent map of $f_1$ is $df - df_0=0$; this proves (3).

  Note that (2) implies that $S_0 \cdot S_1$ is a subgroup.  Since $S$
  is geometrically \gcr, we may apply Proposition
  \ref{prop:twisted-product-gcr} to see that $S_0 \cdot S_1$ is
  geometrically \gcr.  Since $S_1 \lhd S_0 \cdot S_1$ is a normal
  subgroup, it follows from the result of B. Martin (Theorem
  \ref{theorem:martin}) that $S_1$ is \gcr; this proves
  (4).
\end{proof}

\begin{cor}
  \label{cor:implies-main}
  Let $H$ be quasisimple and suppose that the homomorphism $f:H \to G$
  has geometrically \gcr\ image. Then there are uniquely determined
  commuting $\C_G$-homomorphisms $h_0,h_1,\dots,h_r$ and uniquely
  determined non-negative integers $n_0 < n_1 \cdots < n_r$ such that
  $f$ is the twisted-product homomorphism determined by $(\vec h,\vec
  n)$.
\end{cor}

\begin{proof}
  We may use Proposition \ref{prop:frobenius} to find a homomorphism
  $h:H \to G$ and an integer $t \ge 0$ such that $f = h \circ F^t$
  where $F$ is the Frobenius endomorphism of $H$. Moreover, $dh \ne
  0$. If the conclusion of the Theorem holds for $h$, we claim that it
  holds for $f$ as well. Indeed, if $h$ is the twisted-product
  homomorphism determined by the commuting $\C_G$-homomorphisms $\vec
  h = (h_0,h_1,\dots,h_r)$ and the non-negative integers $\vec n = (0
  = n_0 < \cdots < n_r)$, then $f$ is the commuting-product
  homomorphism determined by $\vec h$ and the non-negative integers
  $\vec m = (t < n_1 + t < \cdots < n_r + t)$.  If $f$ had a second
  representation as a commuting-product homomorphism determined by
  $\vec h' = (h_0',\dots,h_t')$ and $\vec m' = (m_0' < \dots < m_t')$
  then using Lemma \ref{lem:untwisted-product} one deduces that $m_0'
  = t$ and one finds a representation of $h$ as the twisted product
  homomorphism determined by $\vec h'$ and $(0 < m_1'-t < \cdots <
  m_t')$.  Thus $\vec h = \vec h'$ and $\vec m = \vec m'$; this proves
  the claim. So we may and will suppose that $df \ne 0$.

  Let us first prove the uniqueness assertion; namely, suppose that
  $\vec h = (h_0,h_1,\dots,h_t)$ and $\vec h' = (h_0',h_1',\dots,h_s)$
  are commuting homomorphisms with $h_i,h_j' \in \C_G$, and suppose
  that $\vec n = (n_0<\dots<n_t)$ and $\vec n' = (n_0'<\dots<n_s')$
  are non-negative integers with $0=n_0 = n'_0$, and suppose that $f =
  \Phi_{\vec h,\vec n} = \Phi_{\vec h',\vec n'}$.  We must argue that
  $s=t$, $\vec h = \vec h'$ and $\vec n =\vec n'$.  We know that
  $df = dh_0 = dh_0'$. Since $h_0 \in \C_G$ is the unique mapping with
  $df = dh_0$ by (C4), we have $h_0=h_0'$.  It then follows that
  \begin{equation*}
    \Phi_{(h_1,\dots,h_t),(n_1<\dots<n_t)} =  
    \Phi_{(h_1',\dots,h_t'),(n_1'<\dots<n_t')}
  \end{equation*}
  so by induction on $\min(s,t)$, we find that $s=t$, $h_i = h_i'$ and
  $n_i = n_i'$ for $1 \le i \le t$; this completes the proof of
  uniqueness.

  For the existence, we choose by (C4) the unique map $f_0 \in \C_G$
  such that $df =df_0.$ We now write $f_1:H \to G$ for the
  homomorphism of Proposition \ref{prop:steinberg-untwist}(1). Thus
  $f$ is given by the rule
  \begin{equation}
    \label{eq:f-rule}
    g \mapsto f_0(g) \cdot f_1(g)
  \end{equation}
  Write $S$ for the image of $f$, and write $S_0$ and $S_1$ for the
  respective images of $f_0$ and $f_1$.

  We proceed by induction on the semisimple rank $r$ of $G$. If $r$ is
  smaller than the rank of the simple group $H$, there are no
  homomorphisms $H \to G$.  If the semisimple rank of $G$ is the same
  as the rank of $H$, then apply Lemma \ref{lem:equal-rank} to $S_0
  \subset G'$, where $G'$ is the derived group of $G$. One deduces
  that $C_{G'}(S_0)$ has no non-trivial torus, hence that any torus in
  $C_G(S_0)$ is central in $G$. Since $\SL_2$ is its own derived
  group, $\image f_1$ lies in $G'$; thus $\image f_1$ is contained in
  $C_{G'}(S_0)$ by Proposition \ref{prop:steinberg-untwist}(2) and it
  follows that the map $f_1$ is trivial. We conclude in this case that
  $f=f_0 \in \C_G$.

  We now suppose that the semisimple rank of $G$ is strictly greater
  than the rank of $H$. Since $S_1 \subset C_G(S_0)$, a maximal torus
  of $S_0$ centralizes $S_1$.  Thus, the image of the \gcr\
  homomorphism $f_1$ lies in some proper Levi subgroup $L$.  Since the
  semisimple rank of $L$ is smaller than that of $G$, we may apply
  induction; we find commuting homomorphisms $h_1,\dots,h_r \in \C_L$,
  and non-negative integers $n_1 < n_2 < \cdots < n_r$ such that $f_1$
  is the twisted-product map determined by $(\vec h,\vec n)$.  Since
  $df_1=0$, we have $0 < n_1$. It follows from (C5) that
  $h_1,\dots,h_r \in \C_G$.

  Since $\image f_0 = S_0 \subset C_G(S_1)$, it follows from
  Proposition \ref{prop:twisted-centralizer} applied to $f_1$ that
  $S_0 \subset C_G(\image h_i)$ for $1 \le i \le r$. Thus
  the homomorphisms $(f_0,h_1,\dots,h_r)$ are commuting.
  In view of \eqref{eq:f-rule}, $f$ is the twisted-product
  homomorphism determined by $(f_0,h_1,\dots,h_r)$ and $(0 < n_1 <
  \cdots < n_r)$.  
\end{proof}

\begin{lem}
  \label{lem:equal-rank}
  Let $X$ and $Y$ be semisimple groups of the same rank, and suppose
  that $X \subset Y$.  Then $C_Y(X)$ contains no non-trivial torus.
\end{lem}

\begin{proof}
  Let $S \subset Y$ be any torus centralizing $X$, and let $T$ be a
  maximal torus of $X$. Since $T$ is centralized by $S$ and is also
  maximal in $Y$, we have $S \subset T$ so that $S \subset X$. Thus
  $S$ is a central torus in $X$. Since $X$ is semisimple, $S$ is
  trivial as required.
\end{proof}

  

We can now prove the following; note that Theorem \ref{theorem:main} is a
special case.
\begin{theorem}
  \label{theorem:real-main}
  Let $G$ be a strongly standard reductive group, and let
  $\Phi:\SL_2 \to G$ be a homomorphism.  If the image of $\Phi$ is
  geometrically \gcr, then there are commuting optimal homomorphisms
  $\vec \phi = (\phi_0,\dots,\phi_r)$ and non-negative integers $\vec
  n = (n_0 < n_1 < \cdots < n_d)$ such that $\Phi$ is the
  twisted-product homomorphism determined by $(\vec \phi,\vec n)$.
  Moreover, $\vec \phi$ and $\vec n$ are uniquely determined by
  $\Phi$.
\end{theorem}

\begin{proof}
  For a strongly standard reductive group $G$, write $\C_G$ for the
  set of optimal homomorphisms $\SL_2 \to G$.  Suppose first that
  $p>2$. Then Theorem \ref{theorem:main} is a consequence of
  Proposition \ref{prop:sl2-conditions-hold} together with Corollary
  \ref{cor:implies-main}.

  Now suppose that $p=2$. Use \eqref{eq:cover-G} to find a central
  isogeny $\pi:G_\sconn \to G$ where the derived group
  of $G_\sconn$ is simply connected. Since $\SL_2$ is simply
  connected, there is a homomorphism $\widetilde{\Phi}:\SL_2 \to
  G_\sconn$ with $\Phi = \pi \circ \widetilde{\Phi}$. It follows from
  Lemma \ref{lem:gcr-isogeny} that $\widetilde{\Phi}$ has
  geometrically \gcr\ image.  
  Proposition \ref{prop:optimal-in-Levi-and-under-isogeny}(2) shows
  that a homomorphism $f:\SL_2 \to G_\sconn$ is optimal if and only if
  $\pi\circ f:\SL_2 \to G$ is optimal.

  If $\widetilde{\Phi}$ is the twisted product homomorphism determined
  by the optimal homomorphisms $\vec \phi = (\phi_0,\dots,\phi_r)$ and
  the non-negative integers $\vec n = (n_0<n_1<\cdots<n_r)$, it is
  then clear that $\Phi$ is the twisted product homomorphism
  determined by the optimal homomorphisms $\vec \phi' = (\pi \circ
  \phi_0,\dots,\pi \circ \phi_r)$ and $\vec n$. Moreover, the
  uniqueness of $\vec \phi$ implies the uniqueness of $\vec \phi'$;
  thus it suffices to prove the theorem after replacing $G$ by
  $G_\sconn$. So we now assume that the derived group of $G$ is simply
  connected.

  Assume first that $K$ is separably closed. Recall that since $p=2$
  is good for the derived group of $G$, each of its simple factors has type
  $A_m$ for some $m$. Since $G$ is split and simply connected, we find
  that $G \simeq T \times \prod_{i=1}^t G_i$ where $T$ is a central
  torus, and $G_i \simeq \SL(V_i)$ for a vector space $V_i$. Write
  $\pi_i:G \to G_i$ for the $i$-th projection, and $\Phi_i = \pi_i
  \circ \Phi$.  Steinberg's tensor product theorem
  \cite{JRAG}*{Cor.  II.3.17} shows that  $\Phi_i$ may be
  written as a twisted-product homomorphism for a unique collection of
  commuting optimal homomorphisms and a unique increasing list of
  non-negative integers; see \cite{liebeck-seitz}*{Lemma 4.1}. The
  same then clearly holds for $\Phi$, so the Theorem is proved in this
  case.

  For general $K$, the above argument represents the base-changed
  morphism $\Phi_{/\Ksep}$ as the twisted product homomorphism
  $\Phi_{\vec \phi,\vec n}$ for unique commuting optimal
  $\Ksep$-homomorphisms $\vec \phi = (\phi_0,\dots,\phi_t)$ and unique
  $\vec n = (n_0 < n_1 < \dots < n_t)$. Since $\vec \phi$ and $\vec n$
  are unique, we may apply Galois descent to see that each $\phi_i$
  arises by base change from an optimal $K$-homomorphism, and the
  proof is complete.
\end{proof}

\section{Proof of a partial converse to the main theorem}
\label{sec:converse-proof}

In this section, we will prove Theorem \ref{theorem:converse}, which
is a geometric statement -- it depends only on $G$ and $H$ over an
algebraically closed field. Thus we will suppose in this section that
$K$ is algebraically closed, and we write ``\gcr'' rather than
``geometrically \gcr''.

We begin with a result on \gcr\ subgroups.
\begin{prop}
  \label{prop:commuting-gcr}
  Let $G$ be reductive, let $h(G)$ be the maximum Coxeter number of a
  simple quotient of $G$, and suppose that $p>2h(G) -2$.  Let
  $A,B \subset G$ be smooth, connected, and \gcr, and suppose that $B \subset
  C_G(A)$. Then $A \cdot B$ is  \gcr.
\end{prop}

\begin{proof}
  Under our assumptions on $p$, it follows from
  \cite{serre-sem-bourb}*{Corollaire 5.5} that a subgroup $\Gamma \subset G$
  is \gcr\ if and only if the representation of $\Gamma$ on $\Lie(G)$ is
  semisimple.

  Since a smooth, connected \gcr\ subgroup is reductive
  \cite{serre-sem-bourb}*{Prop. 4.1}, the proposition is now a
  consequence of the lemma which follows.
\end{proof}

\begin{lem}
  \label{lem:ss}
  Let $G_1,G_2 \subset \GL(V)$ be connected and reductive, and suppose
  $G_2 \subset C_{\GL(V)}(G_1)$. Then $V$ is semisimple for $G_1 \cdot
  G_2$.
\end{lem}

\begin{proof}
  Write $H=G_1 \cdot G_2$. Since $H$ is a quotient of the reductive
  group $G_1 \times G_2$ by a central subgroup, $H$ is reductive.

  Since $G_1$ and $G_2$ commute, $G_2$ leaves stable the
  $G_1$-isotypic components of $V$. Thus we may write $V$ as a direct
  sum of $H$-submodules which are isotypic for both $G_1$ and $G_2$.
  Thus we may as well assume that $V$ itself is isotypic for $G_1$ and
  for $G_2$.

  Let $B_i \subset G_i$ be Borel subgroups and let $T_i \subset B_i$ be maximal
  tori for $i=1,2$. Note that the choice of a Borel subgroup
  determines a system of positive roots in each $X^*(T_i)$; the
  weights of $T_i$ on $U_i = R_u(B_i)$ are \emph{positive.} Our
  hypothesis means that there are dominant weights $\lambda_i \in
  X^*(T_i)$ such that each simple $G_i$-submodule of $V$ is isomorphic
  to $L_{G_i}(\lambda_i)$, the simple $G_i$-module with highest weight
  $\lambda_i$.

  Now, $B=B_1\cdot B_2$ is a Borel subgroup of $H$, and $T=T_1\cdot
  T_2$ is a maximal torus of $B$.  Since $T_1 \cap T_2$ lies in the
  center of $H$, one knows that there is a unique character $\lambda
  \in X^*(T)$ such that $\lambda_{\mid T_i} = \lambda_i$ for $i=1,2$.
  Moreover, it is clear that $\lambda$ is dominant.  Put $U=U_1 \cdot
  U_2 = R_u(B)$.

  It follows from \cite{JRAG}*{II.2.12(1)} that there are no
  non-trivial self-extensions of simple $H$-modules; thus the Lemma will follow if
  we show that all simple $H$-submodules of $V$ are isomorphic to
  $L_H(\lambda)$.

  Let $L \subset V$ be a simple $H$-submodule; we claim that $L \simeq
  L_H(\lambda)$. Since $L$ is simple, the fixed point space of $U$ on
  $L$ satisfies $\dim_K L^U =1$ and our claim will follow once we show
  that $L^U \subset L_{T;\lambda}$ since then $L^U=L_{T;\lambda}$
  and $L \simeq L_H(\lambda)$; for all this, see \cite{JRAG}*{Prop.
    II.2.4} [we are writing $L_{T;\lambda}$ for the $\lambda$ weight
  space of the torus $T$ on $L$]. Since $L$ is semisimple and
  $G_1$-isotypic, $L^{U_1} = L_{T_1;\lambda_1}.$ Since $G_2 \subset
  C_{\GL(V)}(G_1)$, $L^{U_1}$ is a $G_2$-submodule.  Since $L^{U_1}$
  is semisimple and isotypic as $G_2$-module, we know that $L^U =
  (L^{U_1})^{U_2} = (L^{U_1})_{T_2;\lambda_2}.$ Thus $L^U \subset
  L_{T_1;\lambda_1} \cap L_{T_2;\lambda_2}$ so indeed $L^U=
  L_{T;\lambda}$ as required.
\end{proof}

Suppose that $H$ is a simple group, and that for each strongly
standard group $G$, one has a set $\C_G$ of homomorphisms $H \to G$
satisfying (C1)--(C5) of \S \ref{sub:setting}.

\begin{theorem}
  \label{theorem:general-converse}
  Let $G$ be strongly standard and assume that $p > 2h(G) - 2$. Let
  $h_0,\dots,h_r$ be commuting $\C_G$-homomorphisms, and let $n_0 < n_1
  < \cdots < n_r$ be non-negative integers. Then the image of the
  twisted-product homomorphism $h$ determined by $(\vec h,\vec n)$ is
  geometrically \gcr.
\end{theorem}

\begin{proof}
  Write $S_i$ for the image of $h_i$, $0 \le i \le r$. By (C1), $S_i$
  is \gcr\ for $0 \le i \le r$. In view of our assumption on $p$, it
  follows from Proposition \ref{prop:commuting-gcr} that the subgroup
  $A=S_0 \cdot S_1 \cdots S_r$ is \gcr.

  Write $X$ for the building of $G$. If $S = \image h$, Corollary
  \ref{cor:twisted-product-in-P} shows that $X^{S} = X^{A}$. Since $A$
  is \gcr, $X^A = X^S$ is not contractible, so that $S$ is \gcr\
  \cite{serre-sem-bourb}*{Th\'eor\`eme 2.1}.
\end{proof}

\begin{theorem}
  \label{theorem:real-converse}
  Let $G$ be a strongly standard reductive group,   suppose that $p >
  2h(G) -2$, let $\vec \phi = (\phi_0,\dots,\phi_d)$ be commuting
  optimal homomorphisms $\SL_2 \to G$, and let $\vec n = (n_0 < n_1 <
  \cdots < n_d)$ be non-negative integers. Then the image of the
  twisted-product homomorphism $\Phi:\SL_2 \to G$ determined by $(\vec
  \phi,\vec n)$ is geometrically \gcr.
\end{theorem}

\begin{proof}
  As in the proof of Theorem \ref{theorem:main}, write $\C_G$ for the
  set of optimal homomorphisms $\SL_2 \to G$ for a strongly standard
  group $G$.  Note that the condition $p>2h(G) - 2$ implies that $p >
  2$. Then Theorem \ref{theorem:converse} is a consequence of
  Proposition \ref{prop:sl2-conditions-hold} and Theorem
  \ref{theorem:general-converse}.
\end{proof}

Of course, Theorem \ref{theorem:converse} is a special case
of the previous result.

\begin{rem}
  \label{rem:converse-classical}
  Let $G$ be one of the following groups: (i) $\GL(V)$,
  (ii) the symplectic group $\SP(V)$, (iii) the orthogonal group
  $\SO(V)$, or (iv) a group of type $G_2$. In cases (ii), (iii) assume
  $p > 2$ while in case (iv) assume that $p>3$; then $p$ is very good
  for $G$. In case (iv), write $V$ for the 7 dimensional irreducible
  module for $G$; thus in each case $V$ is the ``natural'' module
  for $G$. Then a closed subgroup $H \subset G$ is  \gcr\
  if and only if $V$ is semisimple as an $H$-module; see
  \cite{serre-sem-bourb}*{3.2.2}.  Thus, the conclusion of Theorem
  \ref{theorem:converse} holds for $G$ (with no further prime
  restrictions). Indeed, in view of Lemma \ref{lem:ss}, one finds that
  the conclusion of Proposition \ref{prop:commuting-gcr} is valid with
  no further assumption on $p$ by using $V$ rather than the adjoint
  representation of $G$. Now argue as in the proof of Theorem
  \ref{theorem:general-converse} when $p>2$, or just use Steinberg's
  tensor product theorem when $p=2$ (since we are supposing $G=\GL(V)$
  in that case).
\end{rem}

\newcommand\mylabel[1]{#1\hfil}

\end{document}